\documentclass[12pt,a4paper,dvipdfmx]{amsart}

\usepackage{amsmath,amsfonts,amsthm,amssymb,amscd,extarrows}

\usepackage{graphicx}
\usepackage{here}   
\usepackage[stable]{footmisc}
\usepackage{cancel}
\usepackage{ascmac,fancybox}
\usepackage{cases}   
\usepackage{ulem}   
\usepackage[all]{xy}
\usepackage{url}
\usepackage{color}

\usepackage{tikz}
\usetikzlibrary{intersections,calc,arrows.meta,shadows}

\makeatletter
  
  \@addtoreset{equation}{section}
\makeatother

\newtheorem{Theorem}{Theorem}[section]
\newtheorem{Proposition}[Theorem]{Proposition}
\newtheorem{Lemma}[Theorem]{Lemma}

\newtheorem{Remark}[Theorem]{Remark}

\newtheorem{Definition}[Theorem]{Definition}

\def\RMN#1{\uppercase\expandafter{\romannumeral#1}}

\newcommand{\spec}{\mathop{\rm Spec}\nolimits}
\newcommand{\proj}{\mathop{\rm Proj}\nolimits}

\newcommand{\gp}{\mathfrak{p}}
\newcommand{\gq}{\mathfrak{q}}
\newcommand{\gm}{\mathfrak{m}}

\newcommand{\oo}{\mathcal{O}}

\newcommand{\bZ}{\mathbb{Z}}
\newcommand{\bR}{\mathbb{R}}

\newcommand{\bF}{\mathbb{F}}

\newcommand{\bQ}{\mathbb{Q}}
\newcommand{\bP}{\mathbb{P}}

\begin{document}
\title[Symbolic Rees algebras]{
Symbolic Rees algebras of space monomial primes of degree $5$}
\author{Kazuhiko Kurano}
\dedicatory{Meiji University}
\date{}
\thanks{
This work was supported by JSPS KAKENHI Grant Number 22K03256. \\
}
\maketitle


\begin{abstract}
Let $K$ be a field of characteristic $0$.
Let $\gp_K(5,103,169)$ be the defining ideal of the space monomial curve $\{(t^5,t^{103},t^{169})\mid t\in K\}$.
In this paper we shall prove that the symbolic Rees algebra $R_s(\gp_K(5,103,169))$ is not Noetherian, that is, is not finitely generated over $K$.
\end{abstract}

\section{Introduction}

For a prime ideal $Q$ of a commutative ring $A$, the $n$th symbolic power of $Q$ is defined to be $Q^nA_Q\cap A$, and denoted by $Q^{(n)}$.
We call the subring
\[
A[Qt, Q^{(2)}t^2, Q^{(3)}t^3, \ldots ]
\]
of the polynomial ring $A[t]$ the  {\em symbolic Rees algebra} of $Q$, and denoted by $R_s(Q)$.
This is not necessarily finitely generated as an $A$-algebra.
In the case where $A$ is Noetherian, $R_s(Q)$ is Noetherian 
if and only if $R_s(Q)$ is finitely generated over $A$ as a ring.

In commutative algebra, the study of the finite generation (over the base ring) of symbolic Rees algebras began with Cowsik's problem~\cite{Cowsik}. 
Cowsik's problem asks whether the symbolic Rees algebra of a prime ideal in a polynomial ring over a field is finitely generated. 
Cowsik seems to have considered this question in connection with the set-theoretic complete intersection property of affine curves. 
The finite generation of symbolic Rees algebras is also deeply related to other problems in mathematics.

For example, symbolic Rees algebras often appear as Cox rings of algebraic varieties or as their subrings. 
The finite generation of Cox rings of algebraic varieties, or of their subrings, is a fundamental problem lying at the heart of birational geometry.

In 1956, Nagata~\cite{Nagata} gave a counterexample to Hilbert's fourteenth problem, and in that paper he proposed the following conjecture, now known as Nagata's conjecture:
``Let $n$ be a natural number with $n \ge 10$.
If there exists a plane curve of degree $d$ passing through 
$n$ general points in the complex projective plane, each with multiplicity at least $r$, 
then $d> \sqrt{n}r$.''
Nagata proved this conjecture in the case where 
$n$ is a perfect square, and thereby obtained a counterexample to Hilbert's fourteenth problem. 
The finite generation of the symbolic Rees algebra of the space monomial primes defined below is closely related to Nagata's conjecture (Proposition~5.2 in \cite{CK}).

Let $K$ be a field.
Let $a$, $b$, $c$ be pairwise coprime positive integers.
Let $S=K[x,y,z]$ be a polynomial ring over $K$.
Let $\gp_K(a,b,c)$ be the kernel of the $K$-algebra homomorphism
$S \stackrel{\varphi}{\rightarrow} K[T]$ defined by $\varphi(x)=T^a$, 
 $\varphi(y)=T^b$ and  $\varphi(z)=T^c$.
An ideal of this form is called a space monomial prime ideal.
We sometimes denote $\gp_K(a,b,c)$ simply by $\gp$.

If a prime ideal $Q$ is generated by a regular sequence, 
then the symbolic powers of $Q$ coincide with its ordinary powers. 
Hence, in this case, the symbolic Rees algebra of $Q$ coincides with the ordinary Rees algebra of $Q$, and therefore the symbolic Rees algebra is finitely generated. 
The number of elements of a minimal generating set of a space monomial prime ideal is either 2 or 3 by a result of Herzog~\cite{Her}. 
When the number of minimal generators is 3, a space monomial prime ideal can be regarded as the simplest example of a prime ideal in a polynomial ring over a field that is not generated by a regular sequence. 
For this reason, the symbolic Rees algebras of space monomial primes have been studied extensively.

Huneke~\cite{Hu} gave a simple necessary and sufficient condition for the symbolic Rees algebra of a space monomial prime to be finitely generated. 
Cutkosky~\cite{C} showed that the symbolic Rees algebra of a space monomial prime is the Cox ring of a certain algebraic surface, thereby giving a geometric interpretation of Huneke's criterion. 
Furthermore, using techniques from algebraic geometry and singularity theory, he provided a sufficient condition for the symbolic Rees algebra of a space monomial prime to be finitely generated. 
The first examples of non–finitely generated symbolic Rees algebras were given by Goto–Nishida–Watanabe~\cite{GNW}. 
For example, when the characteristic of the base field is zero, the symbolic Rees algebra corresponding to $(25,29,72)$
is not finitely generated. 
Gonzales–Gonzales–Karu~\cite{GK}, \cite{GAGK} introduced techniques from toric geometry, which made it possible to construct many examples.

In what follows, we call the minimum of $a$, $b$, $c$
 the degree of a space monomial prime $\gp_K(a,b,c)$. 
By a theorem of Cutkosky~\cite{C}, it is known that the symbolic Rees algebra of a space monomial prime is finitely generated when the degree is at most 4, or equal to 6.
Whether the symbolic Rees algebra of a space monomial prime of degree 5 is finitely generated has remained unknown, and this has been one of the intriguing open questions for researchers in this area. 
The main theorem of this paper addresses this problem.

\begin{Theorem}\label{maintheorem}
Let $K$ be a field of characteristic $0$.
Then the symbolic Rees algebra $R_s(\gp_K(5,103,169))$ is not Noetherian, that is, is not finitely generated over $S$.
\end{Theorem}

The author does not know any example of a space monomial prime of degree $8$
whose symbolic Rees algebra is not Noetherian.
In the case where the degree is $5$ or $8$,
if one of the minimal generators of $\gp_K(a,b,c)$ gives a negative curve (defined in Definition~\ref{NC} below),
then the symbolic Rees algebra $R_s(\gp_K(a,b,c))$ is Noetherian~\cite{Ebina}, \cite{Matsu}, \cite{Uchi}.
In the case where $(a,b,c)=(7,15,26)$, $(9,13,29)$, $(10,11,27)$, $(11,21,25)$, $(12,13,17)$, $(13,18,25), \ldots$,
one of the minimal generators of $\gp_K(a,b,c)$ gives a negative curve and
the symbolic Rees algebra $R_s(\gp_K(a,b,c))$ is not Noetherian
if the characteristic of $K$ is $0$.

Finally, we describe the structure of this paper. 
In Section~2, we introduce techniques from toric geometry and develop the general theory of negative curves. 
Using these tools, we show that in the case $(5,103,169)$, 
there exists an element 
$f$ of degree $2065$ in the $7$th symbolic power that defines a negative curve. 
In Section~3, we show that if the symbolic Rees algebra of the space monomial prime 
 $(5,103,169)$ over a field of characteristic zero is finitely generated,
  then there exists an element $g$ of degree $17407$ in the $59$th symbolic power such that $f$ and $g$ satisfy Huneke's criterion. 
  In Section~4, we reduce the problem to the case of characteristic $2$ and show that the homogeneous component of degree $17407$ of the $59$th symbolic power is zero, 
  thereby completing the proof of the theorem. 
  The computations in Section~4 were carried out using a computer, but the results were also verified by hand without the use of a computer.

\section{Negative curve}

Suppose that $a$, $b$, $c$ are pairwise coprime positive integers.
We think that $S=K[x,y,z]$ is a graded ring with $\deg(x)=a$, $\deg(y)=b$ and 
$\deg(z)=c$.
Then $\gp_K(a,b,c)$ is a homogeneous ideal of $S$.

We may suppose 
\[
\gp_K(a,b,c) = (x^s-y^{t_1}z^{u_1}, y^t-z^{u_2}x^{s_2}, z^u-x^{s_3}y^{t_3}) = I_2\left(
\begin{array}{ccc}
x^{s_2} & y^{t_3} & z^{u_1} \\ y^{t_1} & z^{u_2} & x^{s_3} 
\end{array}
\right) ,
\]
where $s=s_2+s_3$, $t=t_1+t_3$ and $u=u_1+u_2$ by Herzog~\cite{Her}.
Let $a'$, $b'$, $c'$ be integers satisfying $a'a+b'b+c'c=1$ and put
\begin{equation}\label{T}
T=x^{a'}y^{b'}z^{c'} .
\end{equation}
Then we have
\[
S[x^{-1}, y^{-1}, z^{-1}] = \left( S[x^{-1}, y^{-1}, z^{-1}]_0 \right)[T^{\pm 1}] .
\]
Putting
\[
v=\frac{z^{u_2}x^{s_2}}{y^t}, \ \ w=\frac{z^u}{x^{s_3}y^{t_3}},
\]
we have
\[
S[x^{-1}, y^{-1}, z^{-1}]_0 = K[x^{\pm 1}, y^{\pm 1}, z^{\pm 1}]_0
=K[v^{\pm 1}, w^{\pm 1}] 
\]
as in the proof of Lemma~3.2 in \cite{KN}.
It is easy to see that $v^{\alpha}w^{\beta}T^{n} \in S$ if and only if 
\[
\left\{
\begin{array}{l}
s_2 \alpha -s_3\beta +a'n \ge 0 \\
-t \alpha -t_3\beta +b'n \ge 0 \\
u_2 \alpha +u\beta +c'n \ge 0 .
\end{array}
\right.
\]
Let $\Delta$ be a triangle defined by
\begin{equation}\label{delta}
\left\{
\begin{array}{l}
s_2 x -s_3y +a'\ \ge 0 \\
-t x -t_3y +b' \ge 0 \\
u_2 x +uy +c' \ge 0 
\end{array}
\right.
\end{equation}
in an $x$-$y$ plane.
We define the Erhart ring of $\Delta$ as 
\begin{equation}\label{erhart}
K[\Delta]:=\bigoplus_{n\ge 0}\left(
\bigoplus_{(\alpha,\beta)\in n\Delta \cap \bZ^2} Kv^\alpha w^\beta
\right)T^n .
\end{equation}
Then $K[\Delta]$ naturally coincides with $S=K[x,y,z]$.
Let $\gq$ be the ideal of $K[v^{\pm 1}, w^{\pm 1}] $ generated by 
$v-1$ and $w-1$.
Then we have
\[
\gp S[x^{-1}, y^{-1}, z^{-1}] = \gq S[x^{-1}, y^{-1}, z^{-1}] .
\]
For $m \ge 0$, we have
\begin{align*}
\gp^{(m)} & = \gp^m S[x^{-1}, y^{-1}, z^{-1}] \cap S 
= \gq^m S[x^{-1}, y^{-1}, z^{-1}] \cap S \\
& = \bigoplus_{n\ge 0}\left( \gq^m \cap 
\bigoplus_{(\alpha,\beta)\in n\Delta \cap \bZ^2} Kv^\alpha w^\beta
\right)T^n .
\end{align*}

The following proposition is well-known.
Here we give an outline of a proof for the reader's convenience.

\begin{Proposition}\label{regseq}
\begin{enumerate}
\item
Let $g_1$, $g_2$ be an $S$-regular sequence with $g_1 \in [\gp_K(a,b,c)^{(r_1)}]_{d_1}$ and 
$g_2\in [\gp_K(a,b,c)^{(r_2)}]_{d_2}$ for some positive integers $r_1$, $r_2$, $d_1$, $d_2$.
Then we have 
\[
\frac{d_1}{r_1} \cdot \frac{d_2}{r_2} \ge abc .
\]
\item
The symbolic Rees algebra $R_s(\gp_K(a,b,c))$ is Noetherian if and only if
there exists an $S$-regular sequence $g_1$, $g_2$ with $g_1 \in [\gp_K(a,b,c)^{(r_1)}]_{d_1}$ and 
$g_2\in [\gp_K(a,b,c)^{(r_2)}]_{d_2}$ for some positive integers $d_1$, $d_2$, $r_1$, $r_2$ such that 
\[
\frac{d_1}{r_1} \cdot \frac{d_2}{r_2} = abc .
\]
\end{enumerate}
\end{Proposition}

\proof
\noindent (1)
Take $h$ such that $g_1$, $g_2$, $h$ is a homogeneous system of parameters of $S$.
Then the Poincar\'e series of $S/(g_1,g_2,h)$ is 
\[
\frac{(1-t^{d_1})(1-t^{d_2})(1-t^{\deg(h)})}{(1-t^{a})(1-t^{b})(1-t^{c})} .
\]
Substituting $1$ for $t$, we obtain
\begin{equation}\label{huneke}
\frac{d_1d_2\deg(h)}{abc} = \ell_S(S/(g_1,g_2,h)) = {\rm e}((h), S/(g_1,g_2)) 
 \ge \ell_{S_\gp}(S_\gp/(g_1,g_2)S_\gp) \cdot {\rm e}((h), S/\gp)
\ge r_1r_2 \deg(h) ,
\end{equation}
where the first inequality follows from the additive formula of multiplicities and
the second one follows from Lemma~2.2 in \cite{KN}.

\noindent (2)
Let $g_1$, $g_2$ be a regular sequence satisfying 
$g_1 \in [\gp^{(r_1)}]_{d_1}$,
$g_2\in [\gp^{(r_2)}]_{d_2}$ and 
$\frac{d_1}{r_1} \cdot \frac{d_2}{r_2} = abc$.
By (\ref{huneke}), we know $\ell_S(S/(g_1,g_2,h)) = r_1r_2 \ell_S(S/(h)+\gp)$.
Then $R_s(\gp_K(a,b,c))$ is Noetherian by Huneke's criterion~\cite{Hu}.

Conversely assume that $R_s(\gp_K(a,b,c))$ is Noetherian.
Put $X=\proj S$.
Let $\pi:Y\rightarrow X$ be the blow-up at the point corresponding to $\gp_K(a,b,c)$.
Let $E$ be the exceptional curve.
Then there exist curves $C$ and $D$ on $Y$ such that $C\neq E$, $D\neq E$ and $C\cap D=\emptyset$ by Cutkosky's criterion~\cite{C}.
Let $g_1$ and $g_2$ be the defining equations of $\pi(C)$ and $\pi(D)$, respectively.
Suppose $g_1 \in [\gp^{(r_1)}]_{d_1} \setminus [\gp^{(r_1+1)}]_{d_1}$ and
$g_2 \in [\gp^{(r_2)}]_{d_2} \setminus [\gp^{(r_2+1)}]_{d_2}$.
Since $\pi(C) \cap \pi(D)$ is not empty and $C\cap D=\emptyset$, we have $\sqrt{(g_1,g_2)}=\gp$.  Therefore the first ``$\ge$'' in (\ref{huneke}) is ``$=$''.
Here $\spec K[v^{\pm 1}, w^{\pm 1}]$ is an affine open subset of $X$.
The blow-up of it at $\gq$ is an open subscheme of $Y$.
The equations of $\pi(C)$ and $\pi(D)$ in $\spec K[v^{\pm 1}, w^{\pm 1}]$ are $g_1':=g_1/T^{\deg(g_1)}$ and $g_2':=g_2/T^{\deg(g_2)}$, respectively,
where $T$ is an element defined in (\ref{T}).
Here we put $B=K[v^{\pm 1}, w^{\pm 1}]_\gq$.
By Lemma~2.2 in \cite{KN}, we obtain $\ell_B(B/(g_1',g_2'))=r_1r_2$
since $C\cap D = \emptyset$.
Since $\ell_{S_\gp}(S_\gp/(g_1,g_2)S_\gp)=\ell_B(B/(g_1',g_2'))$,
the second `` $\ge$'' in (\ref{huneke}) is also ``$=$''.
\qed

\begin{Definition}\label{Hunekecondition}
\begin{rm}
We say that the {\em Huneke's condition} is satisfied if 
there exists an $S$-regular sequence $g_1$, $g_2$ with $g_1 \in [\gp_K(a,b,c)^{(r_1)}]_{d_1}$ and 
$g_2\in [\gp_K(a,b,c)^{(r_2)}]_{d_2}$ for some positive integers $d_1$, $d_2$, $r_1$, $r_2$ such that 
$\frac{d_1}{r_1} \cdot \frac{d_2}{r_2} = abc$.
\end{rm}
\end{Definition}

\begin{Definition}\label{NC}
\begin{rm}
A non-zero homogeneous irreducible polynomial $f$ in $S$ is called a {\em negative curve}
if $f$ is contained in $[\gp_K(a,b,c)^{(r_0)}]_{d_0}$ for some positive integers $r_0$, $d_0$ with
$d_0/r_0 < \sqrt{abc}$.
\end{rm}
\end{Definition}

In the case where the characteristic of $K$ is positive,
Cutkosky~\cite{C} proved that $R_s(\gp_K(a,b,c))$ is Noetherian if there exists a negative curve.
In the case where the characteristic of $K$ is $0$,
$R_s(\gp_K(a,b,c))$ is not necessarily Noetherian even if there exists a negative curve
(Goto-Nishida-Watanabe~\cite{GNW}).

\begin{Remark}\label{NCrem}
\begin{rm}
\begin{enumerate}
\item
Suppose that $f$ is a negative curve, that is, $f$ is an irreducible element in 
$[\gp^{(r_0)}]_{d_0}$ with $d_0/r_0 < \sqrt{abc}$.
Suppose $g \in [\gp_K(a,b,c)^{(r)}]_d$ for some positive integers $r$, $d$ with
$d/r < abc/(d_0/r_0)$.
Then $f$ divides $g$  by Proposition~\ref{regseq} (1).

We know that, if a negative curve exists, it is determined uniquely up to multiplication by 
an element of $K^\times$.
\item
Assume that $[\gp_K(a,b,c)^{(r)}]_d\neq 0$ for some positive integers $r$, $d$ with
$d/r < \sqrt{abc}$.
Suppose $0\neq f \in [\gp_K(a,b,c)^{(r)}]_d$.
Consider the irreducible decomposition $f=f_1f_2 \cdots f_n$.
Since the associated graded ring
\begin{equation}\label{gr}
S/\gp \oplus \gp/\gp^{(2)} \oplus \gp^{(2)}/\gp^{(3)} \oplus \cdots
\end{equation}
is an integral domain, one of $f_i$'s is a negative curve.
Therefore there exists a negative curve if and only if $[\gp_K(a,b,c)^{(r)}]_d\neq 0$ for some positive integers $r$, $d$ with $d/r < \sqrt{abc}$.
\item
Let $L/K$ be a field extension.
It is easy to see 
\[
[\gp_K(a,b,c)^{(r)}]_d\otimes_KL = [\gp_L(a,b,c)^{(r)}]_d
\]
for any $r$, $d$.

We know that $R_s(\gp_K(a,b,c))$ is Noetherian if and only if
so is $R_s(\gp_L(a,b,c))$.
In particular, the finite generation of $R_s(\gp_K(a,b,c))$ depends only on ${\rm ch}(K)$, $a$, $b$, $c$.

There exists a negative curve in the case $K$, $a$, $b$, $c$ 
if and only if so does in the case $L$, $a$, $b$, $c$. 
Furthermore assume that $K$ is a prime field.
If $f$ is a negative curve  in the case $K$, $a$, $b$, $c$, then
$f$ is also a negative curve in the case $L$, $a$, $b$, $c$.
If a negative curve exists, we may assume that it is a polynomial over the prime field and
it is absolutely irreducible.
\item
Suppose $f$ is a negative curve in $[\gp_K(a,b,c)^{(r_0)}]_{d_0}$ for some positive integers $r_0$, $d_0$ with $d_0/r_0 < \sqrt{abc}$.
Then $r_0$ is also uniquely determined as follows:
We may assume that the coefficients of $f$ are in the prime field.
Assume that $f \in [\gp_K(a,b,c)^{(r_0+1)}]_{d_0}$ and $d_0/r_0 < \sqrt{abc}$.
It is easy to see that $\frac{\partial f}{\partial x}$, $\frac{\partial f}{\partial y}$, $\frac{\partial f}{\partial z}$ are contained in $\gp_K(a,b,c)^{(r_0)}$.
If $\frac{\partial f}{\partial x}$ is not $0$, $f$ divides it and 
the degree of $\frac{\partial f}{\partial x}$ is less than that of $f$.
It is a contradiction.
Thus we know $\frac{\partial f}{\partial x} = \frac{\partial f}{\partial y} = \frac{\partial f}{\partial z}=0$.
Then we may assume that the characteristic of $K$ is a prime number $p$.
Since the coefficients of $f$ are in the prime field, 
there exists a polynomial $f'$ such that $f=f'^p$.
It is a contradiction since $f$ is irreducible.
\item
Suppose $f$ is a negative curve in $[\gp_K(a,b,c)^{(r_0)}]_{d_0}$ for some positive integers $r_0$, $d_0$ with $d_0/r_0 < \sqrt{abc}$.
Suppose $r>0$ and $d\ge 0$.
If $d/r < d_0/r_0$, then $[\gp^{(r)}]_d=0$ as follows:
Assume $0\neq g \in [\gp^{(r)}]_d$ and $d/r < d_0/r_0$.
Then, by (1), $f$ divides $g$.
Put $g=fh$ and $h \in [\gp^{(r')}]_{d'}\setminus [\gp^{(r'+1)}]_{d'}$.
Here remark that $r>r_0$ since $d\ge d_0$ and $d/r < d_0/r_0$.
Then $g = fh \in [\gp^{(r_0+r')}]_{d_0+d'}\setminus [\gp^{(r_0+r'+1)}]_{d_0+d'}$
since the associated graded ring (\ref{gr}) is an integral domain and $f \not\in \gp^{(r_0+1)}$ by (4).
Then $r \le r_0+r'$ and $d=d_0+d'$.
Theefore
\[
\frac{d'}{r'} \le \frac{d-d_0}{r-r_0} < \frac{d_0}{r_0} .
\]
Thus $g$ is divisible by $f^n$ for any $n>0$.
It is a contradiction.
\item
Suppose that Huneke's condition is satisfied, that is, 
there exists an $S$-regular sequence $g_1$, $g_2$ with $g_1 \in [\gp_K(a,b,c)^{(r_1)}]_{d_1}$ and 
$g_2\in [\gp_K(a,b,c)^{(r_2)}]_{d_2}$ ($d_1$, $d_2$, $r_1$, $r_2$ are positive integers) such that 
\[
\frac{d_1}{r_1} \cdot \frac{d_2}{r_2} = abc .
\]
We may assume 
\[
\frac{d_1}{r_1} \le \sqrt{abc} \le \frac{d_2}{r_2} .
\]

If a negative curve does not exist, then we have
\[
\frac{d_1}{r_1} = \frac{d_2}{r_2} = \sqrt{abc}
\]
by (2).

Assume that there exists a negative curve $f$ in $[\gp_K(a,b,c)^{(r_0)}]_{d_0}$ for some positive integers $r_0$, $d_0$ with $d_0/r_0 < \sqrt{abc}$.
We have 
\[
\frac{d_0}{r_0} \le \frac{d_1}{r_1}
\]
by (5).
If $\frac{d_0}{r_0} < \frac{d_1}{r_1}$, then both of $g_1$ and $g_2$ are divisible by $f$.
It is a contradiction.
Therefore we obtain 
\[
\frac{d_0}{r_0} = \frac{d_1}{r_1} .
\]
In this case, we can prove $g_1=cf^n$ for a positive integer $n$ and $c\in K^\times$.
\item
Let $g_1$, $g_2$ be homogeneous elements of $S$.
Then  $g_1$, $g_2$ satisfies Huneke's condition if and only if there exists positive integers $\ell_1$ and $\ell_2$ such that $g_1^{\ell_1}$, $g_2^{\ell_2}$ satisfies Huneke's condition.

Assume that there exists a negative curve $f$.
If there exist $g_1$, $g_2$ satisfying Huneke's condition, we may assume that $g_1$ is a power of $f$ by (6). 
Then $f$, $g_2$ also satisfies Huneke's condition.
\item
Let $L/K$ be a field extension.
Assume that there exists a negative curve $f$.
For positive integers $r$, $d$,
there exists $g$ in $[\gp_K(a,b,c)^{(r)}]_d$ such that $f$, $g$ satisfies Huneke's condition if and only if there exists $g'$ in $[\gp_L(a,b,c)^{(r)}]_d$ such that $f$, $g'$ satisfies Huneke's condition.
\end{enumerate}
\end{rm}
\end{Remark}

For a subset $Q$ of $\bZ^2$, we put
\[
KQ := \bigoplus_{(\alpha,\beta)\in Q}Kv^\alpha w^\beta 
\subset K[v^{\pm 1}, w^{\pm 1}]  .
\]
For a finite set $M$, ${}^\#(M)$ denotes the number of elements contained in $M$.

The following proposition is essentially proved in \cite{Miya}.
Here we give an outline of a proof for the reader's convenience.

\begin{Proposition}\label{Miyahara}
Let $K$ be a field of characteristic $0$.
Let $Q$ be a subset of $\bZ^2$.
Let $L$ be a line in $\bR^2$ such that $L$ contains infinitely many points of $\bZ^2$.
Put $Q'=Q \setminus (L\cap Q)$.
Then, for each positive integer $m$, there exists a $K$-linear map
\[
\psi_m : KQ \cap \gq^m \longrightarrow KQ' \cap \gq^{m-1}
\]
satisfying the following conditions:
\begin{enumerate}
\item
The map $\psi_m$ is surjective if ${}^\#(L\cap Q) \ge m$.
\item
The map $\psi_m$ is injective if ${}^\#(L\cap Q) \le m$.
\item
For $(\alpha_0,\beta_0)\in Q'$ and $\xi \in KQ\cap \gq^m$,
the coefficient of $v^{\alpha_0}w^{\beta_0}$ in $\xi$ is $0$ if and only if 
that in $\psi_m(\xi)$ is $0$.
\end{enumerate}
\end{Proposition}

\proof
First of all, remember that, for $h \in K[v^{\pm 1}, w^{\pm 1}]$, the following two conditions are equivalent:
\begin{itemize}
\item $h \in \gq^m$.
\item For non-negative integers $p,\ q$ with $p+q <m$,
\begin{equation}\label{eq}
\left. \frac{\partial^{p+q}h}{\partial v^p \partial w^q} \right|_{(v,w)=(1,1)}=0 .
\end{equation}
\end{itemize}

Considering a multiplication of $v^{\alpha'}w^{\beta'}$ and an action of ${\rm GL}(2,\bZ)$, we may assume that $L$ is the $y$-axis.

It is enough to construct a $K$-isomorphism $\psi_m : KQ \cap \gq^m \longrightarrow KQ' \cap \gq^{m-1}$ satisfying the above (3) in the case ${}^\#(L\cap Q) = m$.

Suppose
\[
L\cap Q = \{ (0,\beta_1), (0,\beta_2), \ldots, (0,\beta_m) \} ,
\]
where $\beta_1<\beta_2< \cdots < \beta_m$.
Remark 
\[
Q'=Q\setminus (L\cap Q)
=\{ (\alpha,\beta)\in Q \mid \alpha \neq 0 \} .
\]

If $\xi$ is in $\gq^m$, then $v \frac{\partial \xi}{\partial v}$ is in $\gq^{m-1}$.
We define a $K$-linear map $\psi_m: KQ \cap \gq^m \longrightarrow KQ' \cap \gq^{m-1}$
by $\psi_m(\xi) = v \frac{\partial \xi}{\partial v}$.
It is easy to check that $\psi$ satisfies the above (3).

We shall prove that $\psi_m$ is an isomorphism.

First we shall prove that $\psi_m$ is injective.
Take $\xi \in {\rm{Ker}} \psi_m$.
Then $\xi$ satisfies $\dfrac{\partial \xi}{\partial v}=0$.
Therefore we may assume
$$\xi = \sum_{j=1}^m c_jw^{\beta_j} ,$$
where $c_j\in K$ for each $j$.
Since $\xi\in \gq^m$, 
we know
$$ \frac{\partial^u \xi}{\partial w^u}(1,1)=0$$
for $u=0,1,\ldots,m-1$.
Hence we obtain 
\begin{equation}\label{3.3}
0 = \frac{\partial^u \xi}{\partial w^u}(1,1)=\sum_{j=1}^m c_j \beta_j(\beta_j-1)\cdots(\beta_j-u+1)
\end{equation}
for $u=0,1,\ldots,m-1$.
The matrix corresponding to the simultaneous equations (\ref{3.3})
is the following:
$$
\begin{pmatrix}
 1& 1 & \cdots  &1 \\ 
 \beta_1 & \beta_2 & \cdots  &\beta_m \\
 \beta_1(\beta_1-1) & \beta_2(\beta_2-1) & \cdots  &\beta_m(\beta_m-1) \\
 \vdots & \vdots  & \ddots  & \vdots \\ 
 \beta_1(\beta_1-1)\cdots(\beta_1 -m+2) & \beta_2(\beta_2-1)\cdots(\beta_2 -m+2) & \cdots  &\beta_m(\beta_m-1)\cdots(\beta_m-m+2) \\ 
 \end{pmatrix}
$$
After some elementary transformations of rows, we obtain 
$$ \begin{pmatrix}
   1 & \cdots & 1 \\   \beta_1 & \cdots & \beta_m \\ \vdots & & \vdots \\ \beta^{m-1}_1 & \cdots & \beta^{m-1}_m 
\end{pmatrix} .
$$
Since this matrix is invertible,
 (\ref{3.3}) implies $c_j=0$ for $j=1,2,\ldots, m$.
 
Next we shall prove that  $\psi_m$ is surjective.
Take 
\begin{equation*}
\eta = \sum_{(\alpha,\beta)\in Q' } b_{(\alpha, \beta)}v^\alpha w^\beta \in kQ' \cap \gq^{m-1} ,
\end{equation*}
where $b_{(\alpha, \beta)} \in K$.
We define 
\begin{equation*}
\tilde{\eta}= \sum_{(\alpha,\beta)\in Q'} \frac{1}{\alpha} b_{(\alpha, \beta)}v^\alpha w^\beta \in KQ'.
\end{equation*}
Put
\begin{equation}\label{3.4}
\xi(c_1,c_2,\ldots,c_m) = \tilde{\eta} + \sum_{j=1}^m c_j w^{\beta_j} \in KQ
\end{equation}
for $c_1,c_2, \ldots , c_m \in K$.
Since
\begin{equation} \label{3.5}
\frac{\partial}{\partial v}\xi(c_1,c_2,\ldots,c_m) = \frac{\partial \tilde{\eta}}{\partial v} = v^{-1}\eta \in \gq^{m-1} ,
\end{equation}
we have
$$\eta =v \frac{\partial}{\partial v}\xi(c_1,c_2,\ldots,c_m) .$$
In order to show $\eta \in {\rm Im}(\psi_m)$,
it is sufficient to show $\xi(c_1,c_2,\ldots,c_m)\in \gq^m$ for some $c_1,c_2,\ldots,c_m\in K$.

By (\ref{3.5}) and the equivalence (\ref{eq}), we know 
$$\left. \frac{\partial^{i+j} \xi(c_1,c_2,\ldots,c_m)}{\partial v^i \partial w^j } \right|_{(v,w)=(1,1)}=0$$
for $i$ and $j$ satisfying $i >0$, $j\ge 0$, and $0 \leq i+j <m$.
Therefore it is enough to show 
\[
\dfrac{\partial^u \xi(c_1,c_2,\ldots,c_m)}{\partial w^u} (1,1)= \frac{\partial^u \widetilde{\eta}}{\partial w^u} (1,1) +\sum_{j=1}^m c_j \beta_j(\beta_j-1)\cdots(\beta_j-u+1)=0
\]
for $u=0, 1, \ldots, m-1$ and for some $c_1,c_2,\ldots,c_m\in K$.
That is to say,
it is enough that $c_1,c_2,\ldots,c_m$ satisfies
\begin{align*}
&\begin{pmatrix}
 1& 1 & \cdots  &1 \\ 
 \beta_1 & \beta_2 & \cdots  &\beta_m \\
 \beta_1(\beta_1-1) & \beta_2(\beta_2-1) & \cdots  &\beta_m(\beta_m-1) \\
 \vdots & \vdots  & \ddots  & \vdots \\ 
 \beta_1(\beta_1-1)\cdots(\beta_1 -m+2) & \beta_2(\beta_2-1)\cdots(\beta_2 -m+2) & \cdots  &\beta_m(\beta_m-1)\cdots(\beta_m-m+2) \\ 
 \end{pmatrix}
\begin{pmatrix}
c_1\\
c_2\\
c_3\\
\vdots \\
c_m
\end{pmatrix} \\
&=\  -\begin{pmatrix}
\tilde{\eta}(1,1) \\
\frac{\partial \tilde{\eta}}{\partial w} (1,1) \\
\frac{\partial^2 \tilde{\eta}}{\partial w^2} (1,1)\\
\vdots \\
\frac{\partial^{m-1} \tilde{\eta}}{\partial w^{m-1}} (1,1)
\end{pmatrix} .
\end{align*}
Since the above $m \times m$ matrix is invertible, 
there exist $c_1,c_2,\ldots,c_m$ satisfying the above equation.
\qed

Put $X=\proj S$.
Let $Y$ be the blow-up at the point corresponding to $\gp_K(a,b,c)$.
Assume that $f$ is a negative curve in $[\gp_K(a,b,c)^{(r_0)}]_{d_0}$ with $d_0/r_0 < \sqrt{abc}$.
Let $C$ be the proper transform of $V_+(f)$.
Then $C$ is a curve on $Y$ satisfying
\[
C^2 = \frac{d_0^2}{abc} -r_0^2 <0 .
\]

\begin{Proposition}\label{NC5*103*169}
Let $K$ be a field of characteristic $0$.
\begin{enumerate}
\item
There exists a negative curve $f \in [\gp_K(5,103,169)^{(7)}]_{2065}$,
that is, there exists an irreducible polynomial $f \in [\gp_K(5,103,169)^{(7)}]_{2065}$.
(Here remark $2065/7 < \sqrt{5\cdot 103 \cdot 169}$.)
The coefficient of $x^{413}$ in $f$ is not zero.
\item
Furthermore assume that $K$ is algebraically closed.
Let $C$ be the proper transform of $V_+(f)$
where $f$ is the negative curve as in (1).
Then $C$ is isomorphic to $\bP^1_K$.
\end{enumerate}
\end{Proposition}

\proof
We shall prove (1).
Remark that
\[
\gp_K(5,103,169)=I_2\left(
\begin{array}{lll}
x^{28} & y & z \\ y^2 & z & x^{47}
\end{array}
\right) ,
\]
where $I_2( \ )$ stands for the ideal generated by the $2$-minors of the given matrix.
Consider the triangle $2065\Delta$, where $\Delta$ is the triangle defined as in (\ref{delta}).
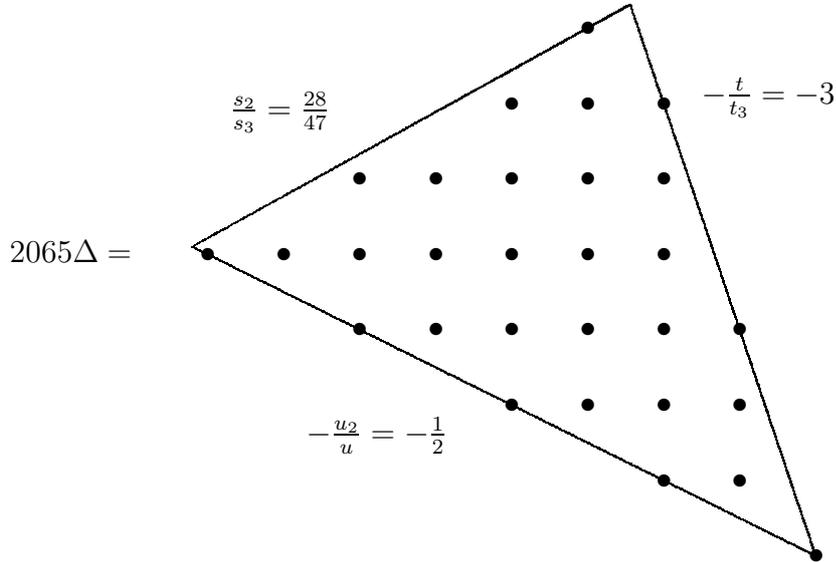
\begin{equation}\label{2065Delta}
{
\setlength\unitlength{1truecm}
  \begin{picture}(11,10)(-2,-5)
\qbezier (-0.2,0.1) (4,-2) (8,-4)
\qbezier (8,-4) (7,-1) (5.56,3.31)
\qbezier (-0.2,0.1) (5,3) (5.56,3.31)
\put(-0.1,-0.1){$\bullet$}
\put(0.9,-0.1){$\bullet$}
\put(1.9,-1.1){$\bullet$}
\put(1.9,-0.1){$\bullet$}
\put(1.9,0.9){$\bullet$}
\put(2.9,-1.1){$\bullet$}
\put(2.9,-0.1){$\bullet$}
\put(2.9,0.9){$\bullet$}
\put(3.9,-2.1){$\bullet$}
\put(3.9,-1.1){$\bullet$}
\put(3.9,-0.1){$\bullet$}
\put(3.9,0.9){$\bullet$}
\put(3.9,1.9){$\bullet$}
\put(4.9,-2.1){$\bullet$}
\put(4.9,-1.1){$\bullet$}
\put(4.9,-0.1){$\bullet$}
\put(4.9,0.9){$\bullet$}
\put(4.9,1.9){$\bullet$}
\put(4.9,2.9){$\bullet$}
\put(5.9,-0.1){$\bullet$}
\put(5.9,1.9){$\bullet$}
\put(5.9,0.9){$\bullet$}
\put(5.9,-3.1){$\bullet$}
\put(5.9,-2.1){$\bullet$}
\put(5.9,-1.1){$\bullet$}
\put(6.9,-1.1){$\bullet$}
\put(6.9,-2.1){$\bullet$}
\put(6.9,-3.1){$\bullet$}
\put(7.9,-4.1){$\bullet$}
  \put(1.3,-2.5){$-\frac{u_2}{u} = -\frac{1}{2}$}
  \put(6.5,2){$-\frac{t}{t_3} = -3$}
  \put(0.3,1.8){$\frac{s_2}{s_3} = \frac{28}{47}$}
    \put(-2.6,-0.1){$2065\Delta =$}
  \end{picture}
}
\end{equation}
Here recall 
\[
S_{2065}=K(2065\Delta \cap \bZ^2)T^{2065}
\]
as in (\ref{erhart}), where $K(2065\Delta \cap \bZ^2)$ is the $K$-vector space spanned by $\{ v^\alpha w^\beta \mid (\alpha,\beta) \in 2065\Delta \cap \bZ^2\}$. 
The bottom lattice point corresponds to $x^{413}$.
The top lattice point corresponds to $y^2z^{11}$.
Here remark 
\[
[\gp_K(5,103,169)^{(7)}]_{2065} = \left( K(2065\Delta \cap \bZ^2) \cap \gq^7
\right) T^{2065} .
\]

Counting from the top, the numbers of lattice points in each row are $1$, $3$, $5$, $7$, $6$, $4$, $2$, $1$.
We know that, by Proposition~\ref{Miyahara},
there exists a polynomial $f$ such that
$$[\gp_K(5,103,169)^{(7)}]_{2065}=Kf ,$$ where the coefficients of $x^{413}$ and $y^2z^{11}$ in $f$ are both non-zero.
(First we set $L$ to be the line containing $7$ lattice points in the fourth row from the top and apply Proposition~\ref{Miyahara}. Next, we set $L$ to be the line containing $6$ lattice points in the fifth row from the top, and so on.)
Hence $f$ is not divided by $x$, $y$ and $z$.
Put $f=\tilde{f}T^{2065}$,
where $$\tilde{f} \in K(2065\Delta \cap \bZ^2) \cap \gq^7 = \bigoplus_{(\alpha,\beta)\in 2065\Delta \cap \bZ^2} Kv^\alpha w^\beta \cap \gq^7$$ as in (\ref{erhart}).
It is enough to show that $\tilde{f}$ is irreducible in $K[v^{\pm 1},w^{\pm 1}]$.

Here consider 
\[
\gp_K(5,11,18)=I_2\left(
\begin{array}{ccc}
x^3 & y & z \\ y^2 & z & x^5
\end{array}
\right) .
\]
Let $\tilde{\Delta}$ be the triangle as in (\ref{delta}) in the case $(a,b,c)=(5,11,18)$.
Then $220\tilde{\Delta}$ and $126\tilde{\Delta}$ are as follows:
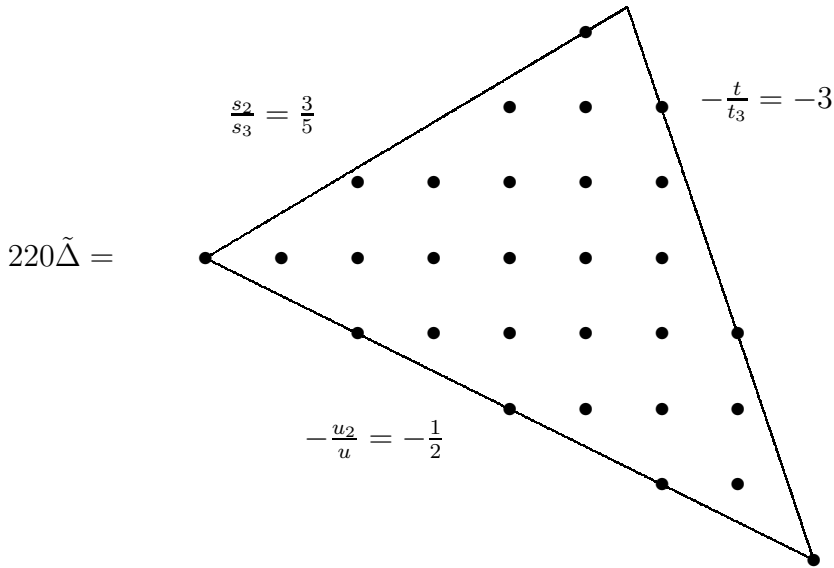
\[
{
\setlength\unitlength{1truecm}
  \begin{picture}(11,10)(-2,-5)
\qbezier (0,0) (4,-2) (8,-4)
\qbezier (8,-4) (7,-1) (5.55,3.33)
\qbezier (0,0) (5,3) (5.55,3.33)
\put(-0.1,-0.1){$\bullet$}
\put(0.9,-0.1){$\bullet$}
\put(1.9,-1.1){$\bullet$}
\put(1.9,-0.1){$\bullet$}
\put(1.9,0.9){$\bullet$}
\put(2.9,-1.1){$\bullet$}
\put(2.9,-0.1){$\bullet$}
\put(2.9,0.9){$\bullet$}
\put(3.9,-2.1){$\bullet$}
\put(3.9,-1.1){$\bullet$}
\put(3.9,-0.1){$\bullet$}
\put(3.9,0.9){$\bullet$}
\put(3.9,1.9){$\bullet$}
\put(4.9,-2.1){$\bullet$}
\put(4.9,-1.1){$\bullet$}
\put(4.9,-0.1){$\bullet$}
\put(4.9,0.9){$\bullet$}
\put(4.9,1.9){$\bullet$}
\put(4.9,2.9){$\bullet$}
\put(5.9,-0.1){$\bullet$}
\put(5.9,1.9){$\bullet$}
\put(5.9,0.9){$\bullet$}
\put(5.9,-3.1){$\bullet$}
\put(5.9,-2.1){$\bullet$}
\put(5.9,-1.1){$\bullet$}
\put(6.9,-1.1){$\bullet$}
\put(6.9,-2.1){$\bullet$}
\put(6.9,-3.1){$\bullet$}
\put(7.9,-4.1){$\bullet$}
  \put(1.3,-2.5){$-\frac{u_2}{u} = -\frac{1}{2}$}
  \put(6.5,2){$-\frac{t}{t_3} = -3$}
  \put(0.3,1.8){$\frac{s_2}{s_3} = \frac{3}{5}$}
  \put(-2.6,-0.1){$220\tilde{\Delta} =$}
  \end{picture}
}
\]
\[
{
\setlength\unitlength{1truecm}
  
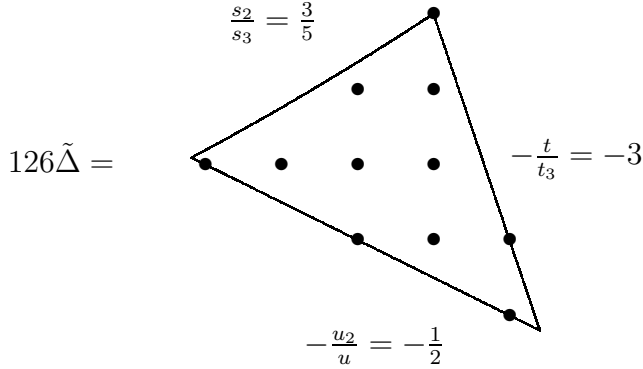
\begin{picture}(11,6)(-2,-3)
\qbezier  (-0.182,0.091) (1.5,1) (3,2)
\qbezier (-0.182,0.091) (2,-1) (4.4,-2.2)
\qbezier (3,2) (4,-1) (4.4,-2.2)
\put(-0.1,-0.1){$\bullet$}
\put(0.9,-0.1){$\bullet$}
\put(1.9,-1.1){$\bullet$}
\put(1.9,-0.1){$\bullet$}
\put(1.9,0.9){$\bullet$}
\put(2.9,1.9){$\bullet$}
\put(2.9,-1.1){$\bullet$}
\put(2.9,-0.1){$\bullet$}
\put(2.9,0.9){$\bullet$}
\put(3.9,-2.1){$\bullet$}
\put(3.9,-1.1){$\bullet$}
  \put(1.3,-2.5){$-\frac{u_2}{u} = -\frac{1}{2}$}
  \put(4,0){$-\frac{t}{t_3} = -3$}
  \put(0.3,1.8){$\frac{s_2}{s_3} = \frac{3}{5}$}
  \put(-2.6,-0.1){$126\tilde{\Delta} =$}
  \end{picture}
}
\]
The bottom lattice point in $220\tilde{\Delta}$ corresponds to $x^{44}$.
The bottom (resp.\ top, left) lattice point in $126\tilde{\Delta}$ corresponds to
$xy^{11}$ (resp.\ $z^7$, $x^{23}y$).
We know that, by Proposition~\ref{Miyahara}, there exists a polynomial $g_1$ such that
$$[\gp_K(5,11,18)^{(7)}]_{220}=Kg_1 , $$ where the coefficient of $x^{44}$ in $g_1$ is non-zero.
Put $g_1=\tilde{g_1}T^{220}$,
where $$\tilde{g_1} \in \bigoplus_{(\alpha,\beta)\in 220\tilde{\Delta} \cap \bZ^2} Kv^\alpha w^\beta \cap \gq^7$$ as in (\ref{erhart}).
Then there exists a unit $\epsilon$ in $K[v^{\pm 1},w^{\pm 1}]$
such that $\tilde{f}=\epsilon\tilde{g_1}$.
By Proposition~\ref{Miyahara}, there exists a polynomial $g_2$ such that
$$[\gp_K(5,11,18)^{(4)}]_{126}=Kg_2 , $$ where the coefficients of $xy^{11}$, $z^{7}$ 
and $x^{23}y$ in $g_2$ are non-zero.
Therefore $g_2$ is not divisible by $x$, $y$, $z$.
Put $g_2=\tilde{g_2}T^{126}$, where $\tilde{g_2} \in K(126\Delta \cap \bZ^2)\cap \gq^4$.
Then $\tilde{g_2}$ is  irreducible by Lemma~2.3 in \cite{GAGK1}
since the coefficients of the top and left lattice points are not $0$.
Hence $g_2$ is irreducible.
Since the coefficient of $x^{44}$ in $g_1$ is non-zero, the sequence $g_1$, $g_2$ is $S$-regular.
Since 
\[
\frac{220}{7} \cdot \frac{126}{4} = 5\cdot 11\cdot 18 ,
\]
$g_1$, $g_2$ satisfies Huneke's condition.
Since 
\[
\frac{220}{7} < \frac{126}{4} ,
\]
$g_1$ is congruence to a power of the negative curve in the case $(a,b,c)=(5,11,18)$
by Remark~\ref{NCrem} (6).
Since $220$ is not a multiple of $7$, 
there does not exist a negative curve with $r_0=1$.
Therefore we know that $g_1$ itself is a negative curve in the case $(a,b,c)=(5,11,18)$.
Hence $g_1$ and $\tilde{g_1}$ are irreducible.
Then $\tilde{f}$ is irreducible.
Therefore $f$  is a negative curve in the case $(a,b,c)=(5,103,169)$.

Next we shall prove (2).
Since $K$ is algebraically closed, it is enough to show $H^1(\oo_C)=0$.
By the exact sequence
\[
0 \longrightarrow \oo_Y(-C) \longrightarrow \oo_Y \longrightarrow \oo_C \longrightarrow 0,
\]
we have an isomorphism
\[
H^1(\oo_C) \simeq H^2(\oo_Y(-C)) 
\]
since $H^1(\oo_Y)=H^2(\oo_Y)=0$ by Leray's spectral sequence.
By the Serre duality,
we obtain
\[
H^2(\oo_Y(-C))=H^0(\oo_Y(C+K_Y)) .
\]
Let $H$ be the pullback of $\oo_X(1)$.
Then $H$ and $E$ are minimal generators of ${\rm Cl}(Y)$.
Then $C$ is linearly equivalent to $d_0H-r_0E=2065H-7E$.
Since $K_Y$ is linearly equivalent to $(-5-103-169)H+E$, 
we have $\oo_Y(C+K_Y) \simeq \oo_Y(1788H-6E)$.
Hence 
\[
H^0(\oo_Y(1788H-6E)) =
K(1788\Delta \cap \bZ^2) \cap \gq^6 .
\]
Here remark that $1788\Delta \cap \bZ^2$ corresponds to 
the lattice points in the interior of $2065\Delta$ in (\ref{2065Delta}).
Counting from the top, the numbers of lattice points in each row in the interior of $2065\Delta$ are $2$, $5$, $6$, $4$, $3$, $1$.
By Proposition~\ref{Miyahara}, we obtain 
\[
K(1788\Delta \cap \bZ^2) \cap \gq^6=0 .
\]
We have completed the proof of Proposition~\ref{NC5*103*169}.
\qed

\section{Huneke's condition}

Let $a$, $b$, $c$ be pairwise coprime positive integers.
Remember $X=\proj S$ and $Y$ is the blow-up of $X$ at $\gp_K(a,b,c)$.

In this section, we shall prove the following theorem:

\begin{Theorem}\label{finitegeneration}
Let $K$ be an algebraically closed field of characteristic $0$.
Let $a$, $b$, $c$ be pairwise coprime positive integers.
Assume that there exists a negative curve $f \in [\gp_K(a,b,c)^{(r_0)}]_{d_0}$ with
$d_0/r_0< \sqrt{abc}$.
Let $d_2$ and $r_2$ be positive integers satisfying the following four conditions
\begin{itemize}
\item[(a)]
$\frac{d_0}{r_0} \cdot \frac{d_2}{r_2} = abc$.
\item[(b)]
$f \not\in (y,z)S$ or $a$ divides $d_2$.
\item[(c)]
$f \not\in (z,x)S$ or $b$ divides $d_2$.
\item[(d)]
$f \not\in (x,y)S$ or $c$ divides $d_2$.
\end{itemize}
Assume the following three conditions:
\begin{enumerate}
\item[(1)]
$R_s(\gp_K(a,b,c))$ is Noetherian.
\item[(2)]
Let $C$ be the proper transform of $V_+(f)$.
Then $C$ is isomorphic to $\bP^1_K$.
\item[(3)]
There exists a positive integer $\ell$ such that
$H^1(\oo_Y(d_2H - r_2E-\ell C))=0$.
\end{enumerate}

Then there exists $g_2 \in [\gp_K(a,b,c)^{(r_2)}]_{d_2}$ such that
$f$, $g_2$ satisfies Huneke's condition.
\end{Theorem}

\begin{Remark}\label{algebraic}
\begin{rm}
We can describe Theorem~\ref{finitegeneration} using only algebraic language.
The above condition (2) is equivalent to $[\gp_K(a,b,c)^{(r_0-1)}]_{d_0-a-b-c}=0$
as in the proof of Proposition~\ref{NC5*103*169} (2).
For the condition (3), we know 
$H^1(\oo_Y(d_2H - r_2E-\ell C)) \neq 0$ if $r_2-\ell r_0 \le -2$, and
\[
H^1(\oo_Y(d_2H - r_2E-\ell C)) = 
\left\{
\begin{array}{ll}
0 & r_2 - \ell r_0 = -1, 0 \\ H^2_{\gm} ( \gp_K(a,b,c)^{(r_2-\ell r_0)} )_{d_2-\ell d_0} 
& r_2-\ell r_0 > 0 ,
\end{array}
\right.
\]
where $\gm=(x,y,z)S$ by Leray's spectral sequence.
\end{rm}
\end{Remark}

The author has two different proofs of the above theorem.
The first proof uses the technique of reduction mod $p$,
and this method is employed in \cite{GNW} and \cite{KN}. 
The second proof examines the transition functions of line bundles; 
this approach is used in \cite{K43}. 
We shall present the second proof below.

Before starting to prove the above theorem,
we need to prove the following easy lemma:

\begin{Lemma}\label{lemmacomplete}
Let $K$ be a field of characteristic $0$.
Let $m$ be a positive integer.
Let $A$ be a $K$-algebra and $I$ be a nilpotent ideal of $A$.

Let $\pi: A \rightarrow A/I$ be the natural surjective ring homomorphism.
Assume that $s \in A$ and $r \in (A/I)^\times$ satisfy $r^m = \pi(s)$.

Then there uniquely exists $u \in A$ such that $\pi(u)=r$ and $u^m=s$.
\end{Lemma}

\proof
It is enough to show it in the case where $I^2=0$.
Take $u'\in A^\times$ such that $\pi(u')=r$.
Suppose $\gamma \in I$.
Then $\pi(u'+\gamma)=r$ is satisfied.
We have
\[
\pi((u'+\gamma)^m)=\pi(u'+\gamma)^m=r^m=\pi(s) .
\]
Therefore there exists $\delta \in I$ such that
\[
\delta = (u'+\gamma)^m-s = ((u')^m-s) + m(u')^{m-1}\gamma .
\]
Here, since $\pi((u')^m-s)=\pi(u')^m -\pi(s)=0$, we know $(u')^m-s \in I$.
If $\gamma = -m^{-1}((u')^{-1})^{m-1}((u')^m-s)$, we have $\delta=0$.
Therefore we obtain $$\left(u'-m^{-1}((u')^{-1})^{m-1}((u')^m-s)\right)^m=s .$$
\qed

\vspace{2mm}

\noindent
{\it Proof of Theorem~\ref{finitegeneration}.}
If $g_1$, $g_2$ satisfy Huneke's condition,
then we may assume that one of $g_1$ and $g_2$ is $f$ by Remark~\ref{NCrem} (7).
Assume that there exists $g \in [\gp_K(a,b,c)^{(r)}]_{d}$ such that
$f$, $g$ satisfies Huneke's condition.
Then we have 
\[
\frac{d_0}{r_0} \cdot \frac{d}{r} = abc .
\]
Therefore we have $d/r=d_2/r_2$ by (a).
Replacing $r$, $d$ by suitable multiples, we may assume that there exists a positive integer $m$ such that $r=mr_2$ and $d=md_2$.

We can take finite number of affine open subsets $U_1$, \ldots, $U_n$ of $Y$ satisfying 
\begin{itemize}
\item
$\{ C\cap U_i \}_{i=1}^n$ is an affine open covering of $C$,
\item
$\oo_Y(d_2H-r_2E)|_{U_i} \simeq \oo_{U_i}$ for $i=1,2,\ldots,n$,
\end{itemize}
for the following reasons:
The singular points of $Y$ come from those of $X=\proj S$.
The singular points of $X$ are contained in $\{ V_+(x,y), V_+(y,z), V_+(z,x) \}$.
If $f \not\in (y,z)S$, then $C$ does not pass through $V_+(y,z)$.
Suppose $f \in (y,z)S$.
Since $a$ divides $d_2$, $\oo_Y(d_2H-r_2E)$ is locally free at $V_+(y,z)$.
Thus $\oo_Y(d_2H-r_2E)$ is locally free near $C$.

The divisorial sheaf $\oo_Y(d_2H-r_2E)$ is a line bundle on $U:=\cup_{i}U_i$.
Let $$\{ t_{ij} \in \Gamma(\oo_Y(U_i\cap U_j))^\times \}_{i,j}$$
be the transition function of $\oo_Y(d_2H-r_2E)|_U$.

Let $D_2$ be the proper transform of $V_+(g)$.
Since $D_2$ is an effective divisor such that $D_2 \sim md_2H-mr_2E$
and $D_2\cap C = \emptyset$, we obtain $\oo_Y(md_2H-mr_2E)|_{\ell C} \simeq \oo_{\ell C}$ for any $\ell >0$.
Here $\ell C$ is the closed subscheme defined by the ideal sheaf $\oo_Y(-\ell C)$.
Therefore there exists $s_i \in \Gamma(\oo_Y(U_i\cap \ell C))^\times$
for $i = 1, 2, \ldots, n$ such that
\[
s_i = (t_{ij}|_{\ell C})^ms_j
\]
in $\Gamma(\oo_Y(U_i\cap U_j \cap \ell C))^\times$
for any $i$, $j$.
Therefore we obtain
\[
s_i|_C = (t_{ij}|_{C})^ms_j|_C
\]
in $\Gamma(\oo_Y(U_i\cap U_j \cap  C))^\times$
for any $i$, $j$.

Since $(d_2H-r_2E).C=0$ and $C \simeq \bP^1_K$,
we have $\oo_Y(d_2H-r_2E)|_C \simeq \oo_C$.
Therefore there exists $r_i \in \Gamma(\oo_Y(U_i\cap C))^\times$
for $i = 1, 2, \ldots, n$ such that
\[
r_i = (t_{ij}|_{C})r_j
\]
in $\Gamma(\oo_Y(U_i\cap U_j \cap  C))^\times$
for any $i$, $j$.
Hence we obtain 
\[
r_i^m = (t_{ij}|_{C})^mr_j^m
\]
in $\Gamma(\oo_Y(U_i\cap U_j \cap  C))^\times$
for any $i$, $j$.
Since $\oo_Y(md_2H-mr_2E)|_C \simeq \oo_C$ and $\Gamma(\oo_C)=K$,
there exists $q \in K^\times$ such that
\[
qr_i^m = s_i|_C
\]
in $\Gamma(\oo_Y(U_i\cap C))^\times$ for $i=1,2,\ldots, n$.
Replacing $r_i$ by $q^{1/m}r_i$, we may assume $q=1$.
Then, by Lemma~\ref{lemmacomplete}, 
there exists $u_i \in \Gamma(\oo_Y(U_i\cap \ell C))^\times$
for $i = 1, 2, \ldots, n$ such that $u_i^m=s_i$ and $u_i|_C=r_i$.
By the uniqueness in Lemma~\ref{lemmacomplete}, 
$u_i=(t_{ij}|_{\ell C})  u_j$ is satisfied for any $i$, $j$.
Thus we know 
\begin{equation}\label{OellC}
\oo_Y(d_2H-r_2E)|_{\ell C} \simeq \oo_{\ell C}
\end{equation}
for any $\ell >0$.

Supppose that $\ell$ is a positive integer satisfying the condition (3) in Theorem~\ref{finitegeneration}.
Consider the exact sequence
\[
0 \longrightarrow \oo_Y(-\ell C) \longrightarrow \oo_Y
\longrightarrow \oo_{\ell C} \longrightarrow 0 .
\]
Since $\oo_Y(d_2H-r_2E)$ is a line bundle near $C$, we obtain the exact sequence
\[
0 \longrightarrow \oo_Y(d_2H-r_2E-\ell C) \longrightarrow \oo_Y(d_2H-r_2E)
\longrightarrow \oo_Y(d_2H-r_2E)|_{\ell C} \longrightarrow 0 .
\]
Then we obtain an exact sequence
\[
H^0(\oo_Y(d_2H-r_2E)) \longrightarrow H^0(\oo_Y(d_2H-r_2E)|_{\ell C})
\longrightarrow H^1(\oo_Y(d_2H-r_2E-\ell C)) .
\]
By the condition (3) in Theorem~\ref{finitegeneration},
we know $H^1(\oo_Y(d_2H-r_2E-\ell C))=0$.
Since 
\[
H^0(\oo_Y(d_2H-r_2E)) \twoheadrightarrow
H^0(\oo_Y(d_2H-r_2E)|_{\ell C}) = H^0(\oo_{\ell C})
\twoheadrightarrow H^0(\oo_{C}) =K
\]
by (\ref{OellC}),
there exists a section in $H^0(\oo_Y(d_2H-r_2E))$ which does not vanish at each point of $C$.
Therefore there exists an effective divisor $D$ such that
$D\sim d_2H-r_2E$ and $C\cap D =\emptyset$.
The defining equation of $D$ satisfies our requirement.
\qed

\section{The proof of infinite generation}

In this section, we shall prove Theorem~\ref{maintheorem}.

We shall prove Theorem~\ref{maintheorem} by contradiction.
Assume that $R_s(\gp_K(5,103,169))$ is Noetherian, where $K$ is an algebraically closed field of characteristic $0$.

There exists a negative curve $f$ as in Proposition~\ref{NC5*103*169} (1)
with $f \not\in (y,z)S$.
Therefore $d_2 = 103\cdot169$ and $r_2=59$ satisfy the conditions (a), (b), (c), (d) in Theorem~\ref{finitegeneration}.
The condition (2) in Theorem~\ref{finitegeneration} is satisfied by Proposition~\ref{NC5*103*169} (2).
We shall prove that the condition (3) in Theorem~\ref{finitegeneration} is satisfied
with $\ell=8$.
Here we have 
\[
d_2H-r_2E-\ell C=103\cdot 169H-59E-8( 2065H - 7E)
= 887H-3E .
\]
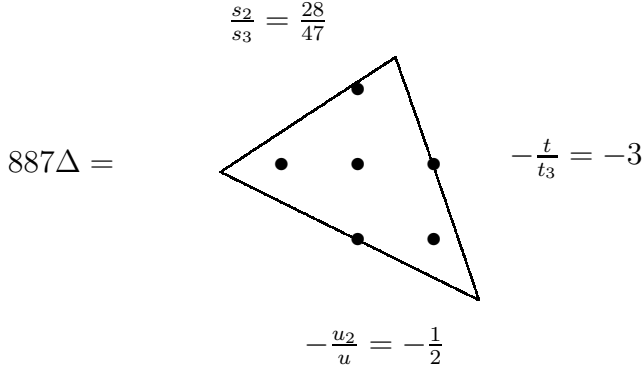
\[
{
\setlength\unitlength{1truecm}
  \begin{picture}(11,6)(-2,-3)
\qbezier  (0.2,-0.1) (1.7,0.9) (2.5,1.42)
\qbezier (0.2,-0.1) (2,-1) (3.6,-1.8)
\qbezier (2.5,1.42) (3,0) (3.6,-1.8)
\put(0.9,-0.1){$\bullet$}
\put(1.9,-1.1){$\bullet$}
\put(1.9,-0.1){$\bullet$}
\put(1.9,0.9){$\bullet$}
\put(2.9,-1.1){$\bullet$}
\put(2.9,-0.1){$\bullet$}
  \put(1.3,-2.5){$-\frac{u_2}{u} = -\frac{1}{2}$}
  \put(4,0){$-\frac{t}{t_3} = -3$}
  \put(0.3,1.8){$\frac{s_2}{s_3} = \frac{28}{47}$}
  \put(-2.6,-0.1){$887\Delta =$}
  \end{picture}
}
\]
We know $\dim_KS_{887}=6$ and $\dim_K[\gp_K(5,103,169)^{(3)}]_{887}=0$
by Proposition~\ref{Miyahara}.
Therefore we obtain $H^1(\oo_Y(887H-3E))=0$.\footnote{
For pairwise coprime positive integers $a$, $b$, $c$, a field $K$,
and positive integers $r$, $d$, 
we have
\[
\frac{r(r+1)}{2}=\dim_KS_d - \dim_KH^0(\oo_Y(dH-rE)) + \dim_KH^1(\oo_Y(dH-rE)) .
\]
}
Hence, if $R_s(\gp_K(5,103,169))$ is Noetherian where $K$ is an algebraically closed field of characteristic $0$, we know that there exists $g_2 \in [\gp_K(5,103,169)^{(59)}]_{17407}$ such that $f$, $g_2$ satisfies Huneke's condition
by Theorem~\ref{finitegeneration}.

Let $\bF_2$ be the prime field of characteristic $2$.
It is easy to see
\[
\dim_K[\gp_K(5,103,169)^{(59)}]_{17407} 
=\dim_\bQ[\gp_\bQ(5,103,169)^{(59)}]_{17407} 
\le \dim_{\bF_2}[\gp_{\bF_2}(5,103,169)^{(59)}]_{17407} .
\]

Then the infinite generation of $R_s(\gp_K(5,103,169))$ immediately comes from the following lemma:

\begin{Lemma}\label{keisan}
$[\gp_{\bF_2}(5,103,169)^{(59)}]_{17407}=0$.
\end{Lemma}

The calculations in the proof of the above lemma were carried out using a computer; however, the author also verified them by hand without using a computer.

\proof
We assume that the base field is $\bF_2$.

Consider the following homogeneous elements of the symbolic powers of $\gp=\gp_{\bF_2}(5,103,169)$:
\begin{align*}
A01 = & y^3 - x^{28}z  \in [\gp]_{309} \\
B01 = & z^2 - x^{47}y  \in [\gp]_{338} \\
C01 = & y^2 z - x^{75} \in [\gp]_{375} \\
D02 = & (y A01 B01 - C01^2)/x^{28} \\
 = & yz^3 + x(\cdots ) \in [\gp^{(2)}]_{610} \\
D03 = & (C01 D02 - A01 B01^2)/x^{19} \\
 = & y^7z + x(\cdots ) \in [\gp^{(3)}]_{890} \\
D04 = & (B01 D03 - A01^2 D02)/x^9 \\
 = & z^7 + x(\cdots ) \in [\gp ^{(4)}]_{1183} \\
D07 = & (D02^2 D03 - A01^3 D04)/x^9 \\
 = & y^2z^{11} + x(\cdots ) \in [\gp^{(7)}]_{2065} \\
D09 = & (A01^2 D07 - D02 D03 D04)/x \\
 = & y^{26} + x(\cdots ) \in [\gp^{(9)}]_{2678} \\
D15 = & (D03 D04^3 - A01 D07^2)/x^2 \\
 = & y^{43} + x(\cdots ) \in [\gp^{(15)}]_{4429} \\
D25 = & (D02 D03^3 D07^2 - D04^4 D09)/x^2 \\
 = & y^{62}z^6 + x(\cdots ) \in [\gp^{(25)}]_{7400} \\
D29 = & (D02 D04^3 D15 - D03^2 D07^2 D09)/x^2 \\
 = & y^{80}z^2 + x(\cdots ) \in [\gp^{(29)}]_{8578} \\
D32 = & (D03^3 D04^2 D15 - D07 D25)/x \\
 = & y^{82}z^6 + x(\cdots ) \in [\gp^{(32)}]_{9460} \\
D37 = & (D03^5 D07 D15 - D04^2 D29)/x \\
 = & y^{98}z^5 + x(\cdots ) \in [\gp^{(37)}]_{10939}  \\
D41 = & (y^2 A01^2 D04 D07^5 - z B01 D02 D03 D07^5 - 
       x^9y  B01^2 D04 D07^5 - x^{10}y  B01 D04^6 D07 D09 \\ & - 
       x^{14}y  D02^2 D37 - x^{15}y  A01^2 D09 D15^2 - x^{38} D04^8 D09 - 
       x^{39} D02 D04^6 D15 - x^{42} D04 D37)/x^{43} \\
 = & y^{116}z + x(\cdots ) \in [\gp^{(41)}]_{12117}  
\end{align*}
 \begin{align*} 
  D43 = & (y^2 A01^2 D04 D07^4 D09 - z B01 D02 D03 D07^4 D09 - 
       x^9y  B01^2 D04 D07^4 D09  \\ & - x^{10}z  A01 D02 D04 D07^3 D15 - 
       x^{11}y  A01 B01 D02 D07^2 D25 - x^{12}y  A01 B01 D02 D07 D32 - 
       \\ & x^{13}y  B01 D02 D03 D37 - x^{14}y  A01^4 D09 D15^2 - 
       x^{22}y  D02 D04 D37 - x^{27}y  A01 D03 D04 D07^5   \\ & - 
       x^{29}y  D03 D04 D07^3 D15 -  x^{39} D02^2 D07^2 D25- 
       x^{40} D02^2 D07 D32 - x^{41} A01^2  D04 D37  \\ &- x^{42} A01^2 D41 - 
       x^{45}  D02 D04^5 D07^3 - x^{46} D02 D03^2 D07^5  - x^{47} D04^7 D15 -
       \\ & x^{48} D03^2 D04^2 D07^2 D15)/x^{51} \\
 = & y^{115}z^5 + x(\cdots ) \in [\gp^{(43)}]_{12690} \\
D49 = & (y B01 D04^4 D07 D25 - z D03^2 D07^4 D15 - 
       x y  B01 D04^4 D32 - x^4y  A01^2 D03^2 D41 \\ & - 
       x^7y  A01^3 D04 D07^6  - x^9y  A01^2 D04 D07^4 D15 - 
       x^{13}y  D03^2 D43 - x^{16}y  D07^7 \\ & - x^{18} B01^2 D04 D07^4 D15  - 
       x^{19} B01 D02 D07^3 D25 - x^{20}  B01 D02 D07^2 D32 \\ & - 
       x^{21} A01 D02^2 D07 D37 - x^{22} B01 D02 D03 D43  - 
       x^{23} A01^4 D15^3 - x^{28} D04^6 D25 \\ & - x^{29} D04^5 D29 - 
       x^{30} D04^3 D37 - x^{31} D02 D04 D43)/x^{32} \\
 = & y^{134}z^4 + x(\cdots ) \in [\gp^{(49)}]_{14478} \\         
D53 = & (y B01 D04^5 D07 D25 - z D03^2 D04 D07^4 D15 - 
       xy  B01 D04^5 D32 - x^4y  D02^2 D49 -   \\ & x^6y  D02 D04^{11} D07 - 
       x^7y  B01 D03 D07^7 - x^{10}y  D03 D04 D07^3 D25 - 
       x^{13}y  D03^2 D04 D43  \\ & - x^{16}y  D04 D07^7 - 
       x^{18} B01^2 D04^2 D07^4 D15 - x^{19} B01 D02 D04 D07^3 D25  \\ & - 
       x^{20}  B01 D02 D04 D07^2 D32 - x^{21} A01 B01 D07^2 D37 - 
       x^{25} B01 D04^6 D07^4 \\ & - x^{28} D04^7 D25 - x^{29} D04^6 D29 - 
       x^{32} D04 D49)/x^{33} \\
 = & y^{152} + x(\cdots ) \in [\gp^{(53)}]_{15656} \\          
D56 = & (y B01 D04^4 D07^2 D25 - z D03^2 D07^5 D15 - 
      xy  B01 D04^4 D07 D32 -x^4 y  A01 D02 D04 D49  \\ & - 
      x^5y  A01 D02 D53 - x^9y  A01^2 D04 D07^5 D15 - 
      x^{13}y  D03^2 D07 D43 - x^{16}y  D07^8  \\ & - x^{18} B01^2 D04 D07^5 D15 -
       x^{19} B01 D02 D07^4 D25 - x^{20} B01 D02 D07^3 D32 \\ & - 
      x^{21} A01 D02^2 D07^2 D37 - x^{22} B01 D02 D03 D07 D43 - 
      x^{23} A01^3 D04 D49 - x^{24} A01^3 D53 \\ & - x^{28} D04^6 D07 D25 - 
      x^{29} D04^6 D32)/x^{33} \\
 = & y^{154}z^4 + x(\cdots ) \in [\gp^{(56)}]_{16538}
\end{align*}

Let $I$ be a homogeneous ideal contained in $\gp^{(59)}$ generated by 
the following $105$ elements:

\vspace{2mm}

{\tiny
$\begin{array}{lllll}
 B01^3 D04^{14}, & B01 D02 D04^{14}, & D04^{13} D07, & B01 D04^{11} D07^2, & B01^2 D04^9 D07^3, \\
D02 D04^9 D07^3, & B01 D02 D04^7 D07^4, & D04^6 D07^5, & B01 D04^4 D07^6, & B01^2 D04^2 D07^7, \\
D02 D04^2 D07^7, & B01 D02 D07^8, & C01 D02 D07^8, & A01 D02 D07^8, & B01 A01^2 D07^8, \\
D03 D07^8, & A01^3 D07^8, & D02 D04^{12} D09, & D03^2 D04 D07^7, & D04^9 D07^2 D09, \\
 B01 D04^7 D07^3 D09, & D03^3 D04^2 D07^6, & D02 D04^5 D07^4 D09, & B01 D02 D04^3 D07^5 D09, & D04^2 D07^6 D09, \\
 B01 D07^7 D09, & C01 D07^7 D09, & D04^{11} D15, & B01 D04^9 D07 D15, & B01^2 D04^7 D07^2 D15, \\
D02 D04^7 D07^2 D15, & B01 D02 D04^5 D07^3 D15, & D04^4 D07^4 D15, & B01 D04^2 D07^5 D15, & B01^2 D07^6 D15, \\
D02 D07^6 D15, & A01 B01 D07^6 D15, & A01 C01 D07^6 D15, & A01^2 D07^6 D15, & D02 D04^8 D25, \\
 B01 D02 D04^6 D07 D25, & D04^5 D07^2 D25, & B01 D04^3 D07^3 D25, & B01^2 D04 D07^4 D25, & D02 D04 D07^4 D25, \\
A01 B01 D04 D07^4 D25, & D03^3 D07^5 D15, & A01^2 D04 D07^4 D25, & B01 D02 D03 D07^4 D25, & D07^5 D09 D15, \\
D02 D04^7 D29, & B01 D02 D04^6 D32, &  D04^5 D07 D32, & B01 D04^3 D07^2 D32, &  B01^2 D04 D07^3 D32, \\
D02 D04 D07^3 D32, & A01 B01 D04 D07^3 D32, & D04^2 D07^3 D15^2, & A01^2 D04 D07^3 D32, & B01 D02 D03 D07^3 D32, \\
A01 D07^4 D15^2, & D02 D04^5 D37, & B01 D02 D04^3 D07 D37, & D04^2 D07^2 D37, & B01 D07^3 D37, \\
 C01 D07^3 D37, & A01 D07^3 D37, & B01 D03 D04 D07^2 D37, & C01 D03 D04 D07^2 D37 , & A01 D03 D04 D07^2 D37, \\
A01^4 D04 D07^2 D37, & D04^4 D43, & B01 D04^2 D07 D43, & B01^2 D07^2 D43, & D02 D07^2 D43, \\
A01 B01 D07^2 D43, & A01 C01 D07^2 D43, & A01^2 D07^2 D43, & A01 B01 D03 D04 D07 D43, & A01 C01 D03 D04 D07 D43, \\
A01^2 D03 D04 D07 D43, & A01^5 D04 D07 D43, & D07^2 D15^3, & D02 D04^2 D49, &  B01 D02 D07 D49, \\
D03^3 D07 D43, & A01 D02 D07 D49, & A01^2 B01 D07 D49, & D03 D07 D49, & A01^3 D07 D49, \\
A01^2 B01 D03 D04 D49, & D03^2 D04 D49, & A01^3 D03 D04 D49, & B01^2 D04 D53, & D02 D04 D53, \\
 B01 D02 D56, & C01 D02 D56, & A01 D02 D56, & A01^2 B01 D56, & D03 D56, \\
A01^3 D56, & A01^2 B01 D03 D53, & D03^2 D53, & A01^3 D03 D53, & A01^6 D53
\end{array}$
}

Let $L$ be a monomial ideal of $k[y,z]$ as follows:

{
\[
L:=
\left(
\begin{array}{lllll}
z^{104}, & yz^{103}, & y^2z^{102}, & y^4z^{101}, & y^6z^{100}, \\
y^7z^{99}, & y^9z^{98}, & y^{10}z^{97} & y^{12}z^{96}, & y^{14}z^{95}, \\
y^{15}z^{94}, & y^{17}z^{93}, & y^{19}z^{92}, & y^{20}z^{91}, & y^{22}z^{90}, \\
y^{23}z^{89}, & y^{25}z^{88}, & y^{27}z^{87}, & y^{28}z^{86}, & y^{30}z^{85}, \\
y^{32}z^{84}, & y^{33}z^{83}, & y^{35}z^{82}, & y^{37}z^{81}, & y^{38}z^{80}, \\
y^{40}z^{79}, & y^{42}z^{78}, & y^{43}z^{77}, & y^{45}z^{76}, & y^{47}z^{75}, \\
y^{48}z^{74}, & y^{50}z^{73}, & y^{51}z^{72}, & y^{53}z^{71}, & y^{55}z^{70}, \\
y^{56}z^{69}, & y^{58}z^{68}, & y^{60}z^{67}, & y^{61}z^{66}, & y^{63}z^{65}, \\
y^{65}z^{64}, & y^{66}z^{63}, & y^{68}z^{62}, & y^{70}z^{61}, & y^{71}z^{60}, \\
y^{73}z^{59}, & y^{74}z^{58}, & y^{76}z^{57}, & y^{78}z^{56}, & y^{79}z^{55}, \\
y^{81}z^{54}, & y^{83}z^{53}, & y^{84}z^{52}, & y^{86}z^{51}, & y^{88}z^{50}, \\
y^{89}z^{49}, & y^{91}z^{48}, & y^{92}z^{47}, & y^{94}z^{46}, & y^{96}z^{45}, \\
y^{97}z^{44}, & y^{99}z^{43}, & y^{101}z^{42}, & y^{102}z^{41}, & y^{104}z^{40}, \\
y^{106}z^{39}, & y^{107}z^{38}, & y^{109}z^{37}, & y^{111}z^{36}, & y^{112}z^{35}, \\
y^{114}z^{34}, & y^{115}z^{33}, & y^{117}z^{32}, & y^{119}z^{31}, & y^{120}z^{30}, \\
y^{122}z^{29}, & y^{124}z^{28}, & y^{125}z^{27}, & y^{127}z^{26}, & y^{129}z^{25}, \\
y^{130}z^{24}, & y^{132}z^{23}, & y^{133}z^{22}, & y^{135}z^{21}, & y^{137}z^{20}, \\
y^{138}z^{19}, & y^{140}z^{18}, & y^{142}z^{17}, & y^{143}z^{16}, & y^{145}z^{15}, \\
y^{147}z^{14}, & y^{148}z^{13}, & y^{150}z^{12}, & y^{152}z^{11}, & y^{153}z^{10}, \\
y^{155}z^{9}, & y^{157}z^{8}, & y^{158}z^{7}, & y^{160}z^{6}, & y^{161}z^{5}, \\
y^{163}z^{4}, & y^{165}z^{3}, & y^{166}z^{2}, & y^{168}z, & y^{170}
\end{array}
\right)
\]
}

Let $D$ be the $i$,$j$th entry of the $21 \times 5$ matrix of generators of $I$.
If $y^kz^l$ is the $i$,$j$th entry of  the $21 \times 5$ matrix of generators of $L$,
then we have 
\[
D \equiv y^kz^l \mod x .
\]
Thus we have $(L,x)S=I+xS$.

Furthermore we have
\[
\ell(S/ I + xS) =
\ell(k[y,z]/L)= 8850 = \left( \frac{59 \times 60}{2} \right) \times 5 = 
e((x), S/\gp^{(59)}) =\ell(S/\gp^{(59)}+xS).
\]
The proof of the second equality in the above equation requires some calculation.
Then we have $I + xS=\gp^{(59)}+xS$.
Since 
$\gp^{(59)}=I+(xS \cap \gp^{(59)})=I+x\gp^{(59)}$, we obtain $I=\gp^{(59)}$ 
by Nakayama's lemma.

The matrix whose entries are the degrees of the generators of $I$ is the following:
\[
\left(
\begin{array}{lllll}
17576, & 17510, & 17444, & 17481, & 17518, \\
17452, & 17489, & 17423, & 17460, & 17497, \\
17431, & 17468, & 17505, & 17439, & 17476, \\
17410, & 17447, & 17484, & 17418, & 17455, \\
17492, & 17426, & 17463, & 17500, & 17434, \\
17471, & 17508, & 17442, & 17479, & 17516, \\
17450, & 17487, & 17421, & 17458, & 17495, \\
17429, & 17466, & 17503, & 17437, & 17474, \\
17511, & 17445, & 17482, & 17519, & 17453, \\
17490, & 17424, & 17461, & 17498, & 17432, \\
17469, & 17506, & 17440, & 17477, & 17514, \\
17448, & 17485, & 17419, & 17456, & 17493, \\
17427, & 17464, & 17501, & 17435, & 17472, \\
17509, & 17443, & 17480, & 17517, & 17451, \\
17488, & 17422, & 17459, & 17496, & 17430, \\
17467, & 17504, & 17438, & 17475, & 17512, \\
17446, & 17483, & 17417, & 17454, & 17491, \\
17425, & 17462, & 17499, & 17433, & 17470, \\
17507, & 17441, & 17478, & 17515, & 17449, \\
17486, & 17523, & 17457, & 17494, & 17428, \\
17465, & 17502, & 17436, & 17473, & 17510
\end{array}
\right) .
\]
Then we know the minimal degree of generators of $I$ is $17410$,
which is the degree of $D03 D07^8$.

Thus we know
\[
[\gp^{(59)}]_{17407}=0 .
\]
\qed

\begin{Remark}
\begin{rm}
Consider the subring
\[
A:=S[A01 t, B01 t, C01 t, D02t^2, D03t^3,\ldots, D56t^{56}]
\]
of the symbolic Rees algebra
\[
R_s(\gp_{\bF_2}(5,103,169)) = S[\gp t, \gp^{(2)} t^2, \gp^{(3)} t^3, \ldots] .
\]
If $n \le 59$, then we have
\[
A\cap St^n = \gp_{\bF_2}(5,103,169)^{(n)}t^n .
\]
We don't need this equation in this paper.
\end{rm}
\end{Remark}

\section*{Acknowledgement}
The author would like to express their sincere gratitude to the members of the Goto Seminar, especially Yasuhiro Shimoda and Koji Nishida.

\vspace{3mm}

\noindent
\begin{tabular}{l}
Kazuhiko Kurano \\
Department of Mathematics \\
Faculty of Science and Technology \\
Meiji University \\
Higashimita 1-1-1, Tama-ku \\
Kawasaki 214-8571, Japan \\
{\tt kurano@meiji.ac.jp} \\
{\tt http://www.isc.meiji.ac.jp/\~{}kurano}
\end{tabular}

\end{document}